\documentclass{amsart}
\usepackage{latexsym,amssymb,amsthm}
\theoremstyle{theorem}
\newtheorem{thm}{Theorem}
\theoremstyle{definition}

\providecommand*{\eu}{\ensuremath{\mathrm{e}}}

\providecommand*{\dd}{\ensuremath{\mathrm{d}}}
\title{Asymptotic formula for Symmetric Involutions}
\author{Yen-chi Roger Lin}
\address{National Taiwan Normal University}
\email{yclin@math.ntnu.edu.tw}
\date{\today}
\begin{document}
\maketitle

The sequence A000898 in OEIS \cite{oeisA000898}:
\begin{quote}
  1, 2, 6, 20, 76, 312, 1384, 6512, 32400, 168992, \dots
\end{quote}
enumerates the numbers of symmetric involutions $I_{2n}^{rc}$, see
Egge \cite{Egge:2007kh}.  It satisfies an obvious recurrence relation:
\begin{equation}
  \label{eq:sir}
  I_0 = 1, \quad I_1 = 2, \quad I_n = 2 \bigl( I_{n-1} + (n-1) I_{n-2}
  \bigr), \,\, \forall\, n \geq 2.
\end{equation}
Using standard techniques it is easy to convert the recurrence
relation to the exponential generating function:
\begin{equation}
  \label{eq:siegf}
  F(z) = \exp \bigl( z^2 + 2z \bigr).
\end{equation}

Down the webpage, OEIS refers its asymptotic formula to Robinson's
paper \cite{Robinson:1976}, but that formula did not yield
satisfactory results.  The computation for $n = 1000$ shows that
Robinson's formula gives the number $8.480 \times 10^{1442}$, which is
about $3.160 \times 10^{-10}$ times the actual number $2.684 \times
10^{1452}$.  In this short note we will give the correct asymptotic
formula for $I_n$.

The main tool we use here is the theorem of Hayman \cite{Hayman:1956},
in which one finds the definition for {\em admissible\/} functions:
\begin{thm}[Hayman]
  \label{thm:Hayman}
  If $f(z)$ is an admissible entire function, with power series $\sum
  a_n z^n$, then
  \begin{equation}
    \label{eq:hayman}
    a_n \sim \frac{ f(r_n) }{ r_n^n \sqrt{ 2\pi b(r_n) } },
  \end{equation}
  where $r_n (>0)$ and $b$ are defined by
  \begin{equation*}
    a(r_n) = n, \quad a(r) = r \, \frac{\dd}{\dd r} \log f(r), \quad
    b(r) = r \, a'(r).
  \end{equation*}
\end{thm}

Here our function $F(z)$ is admissible, hence Theorem~\ref{thm:Hayman}
is applicable.  We have $a(r) = 2r^2 + 2r$, therefore $r_n$ is found
to be:
\begin{equation}
  \label{eq:rn}
  r_n = \frac{ \sqrt{2n+1} - 1 }{2} = \sqrt{ \frac{n}{2} } -
  \frac12 + \frac{1}{4 \sqrt{2n}} + O(n^{-3/2})
\end{equation}
as $n \to \infty$.

From this formula \eqref{eq:rn} for $r_n$, we compute $f(r_n)$,
$r_n^n$, and $b(r_n)$ as follows:

\begin{equation*}
  \begin{split}
    f(r_n) &= \exp( r_n^2 + 2 r_n ) = \exp \Bigl( \frac{1}{2} ( 2 r_n^2
    + 2 r_n ) + r_n \Bigr) \\
    &= \exp \Bigl( \frac{n}{2} + \sqrt{ \frac{n}{2} } - \frac{1}{2}
  \Bigr) \bigl(1 + \frac{1}{4 \sqrt{2n}} + O(n^{-1}) \bigr).
  \end{split}
\end{equation*}

\begin{equation*}
  \begin{split}
    \log r_n^n &= n \log \Bigl( \sqrt{\frac{n}{2}} - \frac{1}{2} +
    \frac{1}{4 \sqrt{2n}} + O(n^{-3/2}) \Bigr) \\
    &= n \Bigl( \log \sqrt{ \frac{n}{2} } + \log \bigl( 1 -
    \frac{1}{\sqrt{2n}} + \frac{1}{4n} + O(n^{-2}) \bigr) \Bigr) \\
    &= n \log \sqrt{ \frac{n}{2} } - \sqrt{ \frac{n}{2} }
    + \frac{1}{12 \sqrt{2n}} + O( n^{-1} ).
  \end{split}
\end{equation*}

\begin{equation*}
  \begin{split}
    \sqrt{b(r_n)} &= \sqrt{4 r_n^2 + 2 r_n} 
    = \sqrt{2 (2 r_n^2 + 2r_n) - 2 r_n} \\
    &= \sqrt{2n - \sqrt{2n} + 1 + O(n^{-1/2})} \\ 
    &= \sqrt{2n} \bigl( 1 - \frac{1}{2 \sqrt{2n}}
    + O(n^{-1}) \bigr).
  \end{split}
\end{equation*}

Finally we put everything back and make use of the Stirling's formula
to derive from \eqref{eq:hayman} that
\begin{equation}
  \label{eq:in}
  I_n = \frac{ \eu^{\sqrt{2n}} }{ \sqrt{2 \eu} } \Bigl( \frac{ 2n }{
    \eu } \Bigr)^{n/2} \bigl( 1 + \frac{\sqrt{2}}{3\sqrt{n}} + O(n^{-1}) \bigr).
\end{equation}

Denote the asymptotic formula on the right-hand side of \eqref{eq:in}
by $I_n^*$.  The following table compares $I_n$ and $I_n^*$ for $n =
10^2, 10^3, 10^4$, and $10^5$:

\begin{equation*}
  \begin{array}{c|c|c|c|c}
    n & 10^2 & 10^3 & 10^4 & 10^5  \\ \hline
    I_n & 1.3506 \cdot10^{99} & 2.6836 \cdot10^{1452} &
    5.3760 \cdot10^{19394} & 4.276309 \cdot10^{243530} \\ \hline
    I_n^* & 1.3520 \cdot10^{99} & 2.6839 \cdot10^{1452} &
    5.3761 \cdot10^{19394} & 4.276313 \cdot10^{243530}
  \end{array}
\end{equation*}

The order of the errors should be $O(n^{-1})$.

\nocite{Bender:1974uq}

\bibliographystyle{amsplain}

\providecommand{\bysame}{\leavevmode\hbox to3em{\hrulefill}\thinspace}
\providecommand{\MR}{\relax\ifhmode\unskip\space\fi MR }
\providecommand{\MRhref}[2]{%
  \href{http://www.ams.org/mathscinet-getitem?mr=#1}{#2}
}
\providecommand{\href}[2]{#2}

\end{document}